\documentclass[11pt]{article}
\usepackage[margin=0.95in]{geometry}
\usepackage{epsfig,ulem,amssymb,soul}

\usepackage{graphicx} 
\usepackage{booktabs}
\usepackage{epstopdf}
\usepackage{amsfonts,amsmath}
\usepackage{bm}
\usepackage{setspace}
\usepackage{comment}
\usepackage{siunitx}
\usepackage{mathrsfs}
\usepackage{multirow}
\usepackage{subcaption}
\usepackage[hang,flushmargin,bottom]{footmisc}
\usepackage{mathtools}
\usepackage{lipsum}
\usepackage{array}
\usepackage{enumerate}
\usepackage{ifthen}
\usepackage[affil-it]{authblk}
\usepackage{cite}
\usepackage{leftidx}
\usepackage{relsize}
\usepackage{braket}
\usepackage{float}
\usepackage{placeins}
\usepackage[colorlinks=true, linkcolor=blue]{hyperref} % Customize hyperlink appearance
\DeclarePairedDelimiter{\ceil}{\lceil}{\rceil}
\usepackage[autostyle]{csquotes}
\usepackage[dvipsnames]{xcolor}

\providecommand{\keywords}[1]{\textbf{\textit{Keywords---}} #1}
\usepackage{hyperref}

\title{On the one-dimensional SPH approximation of fractional-order operators}
\date{}

\author{Khashayar Ghorbani, Fabio Semperlotti\thanks{To whom correspondence should be addressed. Email: fsemperl@purdue.edu }}
\affil{Ray W. Herrick Laboratories, School of Mechanical Engineering,  Purdue University, West Lafayette, IN 47907, USA}
\setstcolor{red}

\usepackage{lineno}
\begin{document}
\maketitle

%\tableofcontents

\begin{abstract}

This work presents a theoretical formalism and the corresponding numerical techniques to obtain the approximation of fractional-order operators over a 1D domain via the smoothed particle hydrodynamics (SPH) method. The method is presented for both constant- and variable-order operators, in either integral or differential forms. Several numerical examples are presented in order to validate the theory against analytical results and to evaluate the performance of the methodology. This formalism paves the way for the solution of fractional-order continuum mechanics models via the SPH method.

\noindent\keywords{Fractional calculus, Constant-order fractional operators, Variable-order fractional operators, Smoothed particle hydrodynamics}
\end{abstract}

\section{Introduction}
\label{sec: Introduction}

The concept of fractional calculus (FC), that is the study of integrals and derivatives of fractional order, has been known to mathematicians and engineers since the early 19th century \cite{tarasov2011fractional,ding2021applications}. Despite the long history, applications of FC to practical problems have not received a significant impulse until the late 20th century \cite{patnaik2020application}. Early applications were mostly seen in control problems \cite{oustaloup1988fractality,podlubny1999fractional} and viscoelastic material simulations \cite{bagley1983theoretical,bagley1986fractional}. They mostly focused on time fractional derivatives due to the innate capability of these operators to effectively capture memory and hereditary effects. In more recent years, fractional operators have also been applied to spatially varying quantities, therefore recognizing the great potential to capture a variety of effects ranging from nonlocal behavior \cite{di2013mechanically,lazopoulos2006non,carpinteri2014nonlocal,patnaik2020ritz,sidhardh2021analysis}, to multiscale effects \cite{patnaik2020application,ding2023fractional,magin2010fractional,chen2004fractional}, to anomalous and hybrid transport \cite{treeby2010modeling,chen2020multi,gomez2016modeling,buonocore2020scattering}.

These applications naturally involve partial differential equations (PDEs) that are fractional in space and, eventually, in time. Their solution on geometries of practical interest requires the development of accurate and robust numerical methods. In this context, various approaches have been introduced to compute fractional-order operators and to solve certain (mostly linear) fractional-order differential equations. These approaches include finite element (FE) \cite{patnaik2020ritz,ding2024multimesh} and meshless \cite{rajan2024element} methods for fractional-order nonlocal elasticity problems, finite difference (FD) scheme for fractional viscoelastic models \cite{galucio2004finite,cortes2007finite}, as well as meshless methods applied to time and/or space-fractional dispersion \cite{yuan2016advanced,pang2015space}, diffusion \cite{lin2018reproducing,Srivastava2021NumericalSO,lin2020petrov,irfan2023rbf}, and wave \cite{raei2021efficient,cheng2018meshless} equations. While these methods have been very effective in solving problems governed by fractional-order equations, there is still a very limited availability of methods capable of addressing problems involving highly nonlinear and discontinuous fields, such as those emerging in large deformation and fracture mechanics problems. In classical continuum mechanics, this class of problems is typically addressed via meshless methods because they are immune from the computational issues following large mesh distortion. Among the meshless methods, smoothed particle hydrodynamics (SPH) provides a powerful approach that is applicable to a wide range of solid and fluid mechanics problems.

SPH was originally introduced by Lucy \cite{lucy1977numerical}, and independently by Gingold and Monaghan \cite{gingold1977smoothed,monaghan1988introduction} to simulate continuum problems in astrophysics. Since then, it has been successfully extended and applied to various engineering problems, including the simulation of compressible and non-compressible fluid flows \cite{cummins1999sph,monaghan1994simulating,monaghan2005smoothed,shao2003incompressible}, free surface fluid flows \cite{monaghan1994simulating,becker2007weakly}, fluid-structure interaction (FSI)\cite{antoci2007numerical,rafiee2009sph,o2021fluid}, large deformation analysis of geomaterials \cite{bui2008lagrangian} and elastodynamics \cite{gray2001sph}, among others. 
%\KG{There are much more papers that can be mentioned here}
One of the primary advantages of SPH techniques is the ease with which different physics can be incorporated within its general formulation \cite{libersky2005smooth}. The fundamental idea behind SPH is based on the integral representation or, equivalently, the kernel approximation of field variables \cite{liu2003smoothed}. By discretizing the computational domain into a set of particles or points, the integral representation of a field variable at a given particle is transformed into a discrete weighted sum. Using this approach, SPH provides effective formulations capable of approximating both field variables and their gradient (or divergence).

The present study focuses on the development of the mathematical formalism and on the numerical validation of SPH method to numerically approximate fractional-order operators. This work will serve as the foundational step to develop SPH techniques to solve space fractional continuum mechanics models.

To apply SPH to evaluate fractional operators, integrals must be computed in addition to the classical gradient computations. While SPH-based approximations of gradients are well established, a corresponding formulation for integrals has not been fully developed. To the best of the author's knowledge, only one study \cite{kisu2013formulation} has previously explored this direction by proposing a SPH-based approximation of the Riemann-Liouville (RL) fractional constant-order (CO) derivative. The present study significantly extends this topic by developing a comprehensive SPH-based methodology to approximate fractional-order operators. Specifically, this work addresses both constant- and variable-order definitions for various types of operators, including the RL and the Caputo definitions. Nevertheless, the methodology is general and can be employed to calculate any other definition of a fractional operator or, for that matter, of an integral operator. This means that this same formulation would allow the calculation of classical integral formalisms arising in nonlocal continuum mechanics \cite{eringen1972nonlocal}. Finally, a modified form of the SPH-based integration is proposed in order to improve the accuracy of the numerical approximations when dealing with this class of intrinsically nonlocal operators.

The remainder of this paper is organized as follows: \S\ref{sec:FC basics} and \S\ref{sec:SPH basics} present the basic definitions of fractional-order operators and the fundamentals of the SPH approach, respectively. \S\ref{sec:SPH to FC} describes the implementation of the SPH approach for fractional-order operators. Then, in \S\ref{sec:singularity}, the issue of kernel singularity is discussed followed by numerical results in \S\ref{sec: results}.

\section{Basic definitions of fractional-order operators}\label{sec:FC basics}
This section provides a brief overview of the basic terminology and mathematical definitions used to define the fractional-order operators that will be used in the following. It is well known that fractional-order operators have many definitions, which are not all equivalent \cite{oldham1974fractional, podlubny1998fractional, kilbas2006theory}. Given that this study is intended to lay the necessary mathematical foundation to enable SPH simulations based on the fractional formulation of mechanics \cite{patnaik2020geometrically,patnaik2022variable,ding2022multiscale,ding2023fractional,ding2024transversely,patnaik2020generalized,patnaik2021fractional,sumelka2014fractional,sumelka2014thermoelasticity,sumelka2015fractional,drapaca2012fractional}, the definitions are restricted to the more common operators for fractional continuum mechanics, namely RL and Caputo. Also, fractional-order operators can be defined as having constant and variable order (as well as combination of them). Therefore, in the following, we review the main definitions following a classification based on the type of order.

\subsection{CO operators}
The most basic definition of a fractional-order operator accounts for a real and constant order. According to Riemann-Liouville, the left- and right-handed RL fractional integrals of a constant order $\alpha$ are defined as \cite{kilbas2006theory}:

\begin{equation}
    %\prescript{}{a}{\mathrm{I}}_x^{m-\alpha}f(x) = \frac{1}{\Gamma(m-\alpha)} \int_a^x (x-x^\prime)^{m-\alpha-1} f(x^\prime) \mathrm{d}x^\prime
    \prescript{RL}{a}{\mathrm{I}}_x^{\alpha}f(x) = \frac{1}{\Gamma(\alpha)} \int_a^x (x-x^\prime)^{\alpha-1} f(x^\prime) \mathrm{d}x^\prime
    \label{eq:RL_Integral_Left}
\end{equation}

\begin{equation}
    %\prescript{}{x}{\mathrm{I}}_b^{m-\alpha}f(x) = \frac{1}{\Gamma(m-\alpha)} \int_x^b (x^\prime-x)^{m-\alpha-1} f(x^\prime) \mathrm{d}x^\prime
        \prescript{RL}{x}{\mathrm{I}}_b^{\alpha}f(x) = \frac{1}{\Gamma(\alpha)} \int_x^b (x^\prime-x)^{\alpha-1} f(x^\prime) \mathrm{d}x^\prime
    \label{eq:RL_Integral_Right}
\end{equation}

where $\Gamma$ is the Gamma function, $x'$ is a dummy variable of integration, and the bounds satisfy the inequalities $x>a$ and $x<b$.
The left- and right-handed RL fractional derivatives can be obtained by applying an $n$-fold integer-order differentiation to the left- and right-handed RL fractional integrals of order $(n-\alpha)$, as follows \cite{kilbas2006theory}: %Eqs.\ref{eq:RL_Integral_Left} and \ref{eq:RL_Integral_Right}:

\begin{equation}
\begin{split}
    %\prescript{RL}{a}{\mathrm{D}}_x^{\alpha}f(x) = \frac{1}{\Gamma(m-\alpha)} \frac{d^m}{dx^m} \int_a^x (x-x^\prime)^{m-\alpha-1} f(x^\prime) \mathrm{d}x^\prime
    \prescript{RL}{a}{\mathrm{D}}_x^{\alpha}f(x) &= \frac{d^n}{dx^n} \left(\prescript{RL}{a}{\mathrm{I}}_x^{n-\alpha}f(x)\right) \\
    &= \frac{1}{\Gamma(n-\alpha)} \frac{d^n}{dx^n} \left(\int_a^x (x-x^\prime)^{n-\alpha-1} f(x^\prime) \mathrm{d}x^\prime\right)
    \label{eq:RL_Derivative_Left}
\end{split}    
\end{equation}

\begin{equation}
\begin{split}
    %\prescript{RL}{x}{\mathrm{D}}_b^{\alpha}f(x) = \frac{(-1)^m}{\Gamma(m-\alpha)} \frac{d^m}{dx^m} \int_x^b (x^\prime-x)^{m-\alpha-1} f(x^\prime) \mathrm{d}x^\prime
    \prescript{RL}{x}{\mathrm{D}}_b^{\alpha}f(x) &= \frac{d^n}{dx^n} \left(\prescript{RL}{x}{\mathrm{I}}_b^{n-\alpha}f(x)\right) \\
    &= \frac{(-1)^n}{\Gamma(n-\alpha)} \frac{d^n}{dx^n} \left(\int_x^b (x^\prime-x)^{n-\alpha-1} f(x^\prime) \mathrm{d}x^\prime \right)
    \label{eq:RL_Derivative_Right}
\end{split}
\end{equation}

where $\frac{d^n}{dx^n}(\cdot)$ indicates an integer-order derivative of order $n$ with $n=\ceil{\alpha}$.

Similarly, the left- and right-handed Caputo fractional derivatives are obtained by applying the left- and right-handed RL fractional integrals of order $(n-\alpha)$ to the integer-order derivative $\frac{d^n f(x)}{dx^n}$, which results in \cite{kilbas2006theory}:

\begin{equation}
\begin{split}
    %\prescript{C}{a}{\mathrm{D}}_x^{\alpha}f(x) = \frac{1}{\Gamma(m-\alpha)} \int_a^x (x-x^\prime)^{m-\alpha-1} \frac{d^m f(x^\prime)}{dx^m} \mathrm{d}x^\prime
    \prescript{C}{a}{\mathrm{D}}_x^{\alpha}f(x) &=  \prescript{RL}{a}{\mathrm{I}}_x^{n-\alpha}\left(\frac{d^n f(x)}{dx^n}\right) \\
    &= \frac{1}{\Gamma(n-\alpha)} \int_a^x (x-x^\prime)^{n-\alpha-1} \frac{d^n f(x^\prime)}{dx^n} \mathrm{d}x^\prime
    \label{eq:Caputo_CO_Derivative_Left}
\end{split}
\end{equation}

\begin{equation}
\begin{split}
    %\prescript{C}{x}{\mathrm{D}}_b^{\alpha}f(x) = \frac{(-1)^m}{\Gamma(m-\alpha)} \int_x^b (x^\prime-x)^{m-\alpha-1} \frac{d^m f(x^\prime)}{dx^m} \mathrm{d}x^\prime
    \prescript{C}{x}{\mathrm{D}}_b^{\alpha}f(x) &=  \prescript{RL}{x}{\mathrm{I}}_b^{n-\alpha}\left(\frac{d^n f(x)}{dx^n}\right) \\
    &= \frac{(-1)^n}{\Gamma(n-\alpha)} \int_x^b (x^\prime-x)^{n-\alpha-1} \frac{d^n f(x^\prime)}{dx^n} \mathrm{d}x^\prime
    \label{eq:Caputo_CO_Derivative_Right}
\end{split}
\end{equation}

While both left- and right-handed RL and Caputo fractional operators have been defined above for completeness and to provide a general foundation, we focus exclusively on the left-handed forms in the remainder of this paper. Extending the application of the SPH-based formulation to the right-handed operators is just a matter of using the appropriate initial definitions given in Eqs (\ref{eq:RL_Integral_Right}), (\ref{eq:RL_Derivative_Right}), and (\ref{eq:Caputo_CO_Derivative_Right}).

\subsection{Variable-order (VO) operators}\label{sec:VO}
In a seminal work, Samko and Ross \cite{samko1993integration} addressed the concept of generalized FC that leads to the definition of VO operators. In this case, the order can be a function of an independent variable, hence potentially exhibiting a rather elaborate memory behavior. The most commonly used independent variables are time and space, but realistically any other quantity (e.g. temperature, stress, etc.) could be chosen based on the physical behavior of interest. When the independent variable is time, the VO operator retains memory of past events and, more specifically, of the past orders. When the independent variable is space, the VO operator retains information at points other than the current evaluation point and of their spatial evolution, therefore effectively describing a spatial memory of the nonlocal behavior \cite{lorenzo2007initialization,lorenzo2002variable}.

Lorenzo and Hartley \cite{lorenzo2002variable} identified three main types of VO operators that result in a different memory behavior. Accordingly, Type-I VO assumes an order variation of the form $\alpha(x)$, where $x$ can be seen as a generic variable. Because the order function depends only on the instantaneous value of the variable $x$ and not on its previous values, this operator definition does not maintain any memory of previous orders. In Type-II and Type-III, the VO are described in the form $\alpha(x')$ and $\alpha(x,x')$, respectively, which result in a different strength of the memory effect. In this study, we only focus on Type-I VO, primarily for two reasons. First, the analytical non-singular forms to be developed for CO operators are also applicable to Type-I VO (see further details in \S\ref{sec:singularity}). Second, the available analytical solutions for the CO Caputo derivative and CO RL integral, presented in \ref{sec:Appendix A}, can also be used for Type-I VO operator. This is due to the fact that, in both cases, the variable of integration is $x'$ rather than $x$, allowing $\alpha (x)$ to be treated as a constant. However, the RL derivative includes a differentiation with respect to $x$, meaning that $\alpha(x)$ can no longer be treated as a constant, therefore making the CO analytical solutions invalid for the RL derivative in Type-I VO.
According to Lorenzo and Hartley \cite{lorenzo2002variable}, a possible definition of the left-handed Type-I VO RL integral is:

\begin{equation}
    \prescript{RL}{a}{\mathrm{I}}_x^{\alpha(x)}f(x) = \frac{1}{\Gamma(\alpha(x))}\int_a^x (x-x^\prime)^{\alpha(x)-1} f(x^\prime) \mathrm{d}x^\prime
    \label{eq:VORL_Integral}
\end{equation}

Based on the definitions given in Eqs.~\ref{eq:RL_Derivative_Left} and \ref{eq:Caputo_CO_Derivative_Left}, it is possible to define the left-handed Type-I VO RL and Caputo derivatives, respectively, as:

\begin{equation}
    \prescript{RL}{a}{\mathrm{D}}_x^{\alpha(x)}f(x) = \frac{d^n}{dx^n}\left(\frac{1}{\Gamma(n-\alpha(x))}\int_a^x (x-x^\prime)^{n-\alpha(x)-1} f(x^\prime) \mathrm{d}x^\prime\right)
    \label{eq:VORL_Integral}
\end{equation}

and

\begin{equation}
    \prescript{C}{a}{\mathrm{D}}_x^{\alpha(x)}f(x) = \frac{1}{\Gamma(n-\alpha(x))}\int_a^x (x-x^\prime)^{n-\alpha(x)-1} \frac{d^nf(x^\prime)}{dx^n} \mathrm{d}x^\prime
    \label{eq:VOCapouto_Integral}
\end{equation}

\section{Fundamentals of SPH approximation}\label{sec:SPH basics}
The SPH approach to the solution of PDEs relies on two main levels of function approximation: 1) the kernel approximation (also referred to as the integral representation), and 2) the particle approximation. In the following, we briefly review the key elements behind these approximations before specializing them for the computation of fractional operators.

\subsection{Integral representation of functions, derivatives, and integrals}
In SPH, a function $f(\boldsymbol{x})$ is approximated in an integral form as \cite{monaghan1988introduction}:

\begin{equation}
    \langle f(\boldsymbol{x}) \rangle = \int_\Omega f(\boldsymbol{x'}) W(\boldsymbol{x-x'},h) d\boldsymbol{x'}
    \label{eq:SPH_int_approx}
\end{equation}

where $\boldsymbol{x'}$ is a dummy variable of integration, $h$ is the smoothing length, and $W(\boldsymbol{x-x'},h)$ is the interpolating kernel function. The angle brackets $ \langle \cdot \rangle$ denote the approximate value of $f(\boldsymbol{x})$ resulting from the integral representation, and $\Omega$ represents the volume of integration containing $\boldsymbol{x}$. 

A comprehensive discussion of how Eq.~\ref{eq:SPH_int_approx} is capable of approximating the function $f(\boldsymbol{x})$ and of the various interpolating kernel functions $W(\boldsymbol{x-x'},h)$ is available in numerous sources and, for brevity, is not reported again here; the interested reader is remanded, as an example, to \cite{monaghan1988introduction,monaghan2005smoothed,liu2003smoothed}.

In the integral representation, the gradient of a function is given by: 

\begin{equation}
    \langle \nabla f(\boldsymbol{x}) \rangle = - \int_\Omega f(\boldsymbol{x'}) \nabla W(\boldsymbol{x-x'},h) d\boldsymbol{x'}
    \label{eq:SPH_grad_approx}
\end{equation}

which has the advantage that only the gradient of the kernel function is required. %Note that Eq.~\ref{eq:SPH_grad_approx} is appropriate for points inside the support domain, but it does require some attention for points on the boundary.In classical SPH,
The approximations of a function and its gradient are typically sufficient to estimate all necessary terms when solving the governing equations of continuum mechanics using the SPH method. However, when dealing with fractional operators (and, hence, fractional differential equations), an explicit evaluation of the integral is also needed. 

To integrate a function via the SPH method, we start from Eq.~\ref{eq:SPH_int_approx} and perform an integration with respect to $\boldsymbol{x}$ over a generic domain $\Pi$: 

\begin{equation}
    %\int_\Pi \langle f(\boldsymbol{x}) \rangle d\boldsymbol{x} = \int_\Pi \int_\Omega f(\boldsymbol{x'}) W(\boldsymbol{x-x'},h) d\boldsymbol{x'} d\boldsymbol{x}
    \mathcal{I}(f(\boldsymbol{x}))=\int_\Pi f(\boldsymbol{x}) d\boldsymbol{x} = \int_\Pi \int_\Omega f(\boldsymbol{x'}) W(\boldsymbol{x-x'},h) d\boldsymbol{x'} d\boldsymbol{x}
    \label{eq:SPH_integral_kernel}
\end{equation}

By exchanging the order of integration on the right-hand side, the approximated form of the integrated function is obtained as:

\begin{equation}
   %\langle \mathcal{I}(f(\boldsymbol{x})) \rangle = \int_{\Pi^*} \langle f(\boldsymbol{x}) \rangle d\boldsymbol{x} = \int_{\Pi^*} f(\boldsymbol{x'}) \Tilde{W}(\boldsymbol{x-x'},h) d\boldsymbol{x'}
    \langle \mathcal{I}(f(\boldsymbol{x})) \rangle = \int_{\Pi^*} f(\boldsymbol{x'}) \Tilde{W}(\boldsymbol{x-x'},h) d\boldsymbol{x'}
    \label{eq:SPH_integral_kernel_final}
\end{equation}

where $\Tilde{W}(\boldsymbol{x-x'},h)=\int_{\Omega^*} W(\boldsymbol{x-x'},h) d\boldsymbol{x}$. Note that, depending on the nature of the integration domains and on their dependence on the integration variables, exchanging the order of integration may also lead to new integration bounds. Consequently, two new generic domains $\Pi^*$ and $\Omega^*$ are introduced to properly account for this eventuality after reordering.

\subsection{Particle approximation}
\begin{figure}[h!]
    \centering
    \includegraphics[width=0.55\textwidth, keepaspectratio]
    {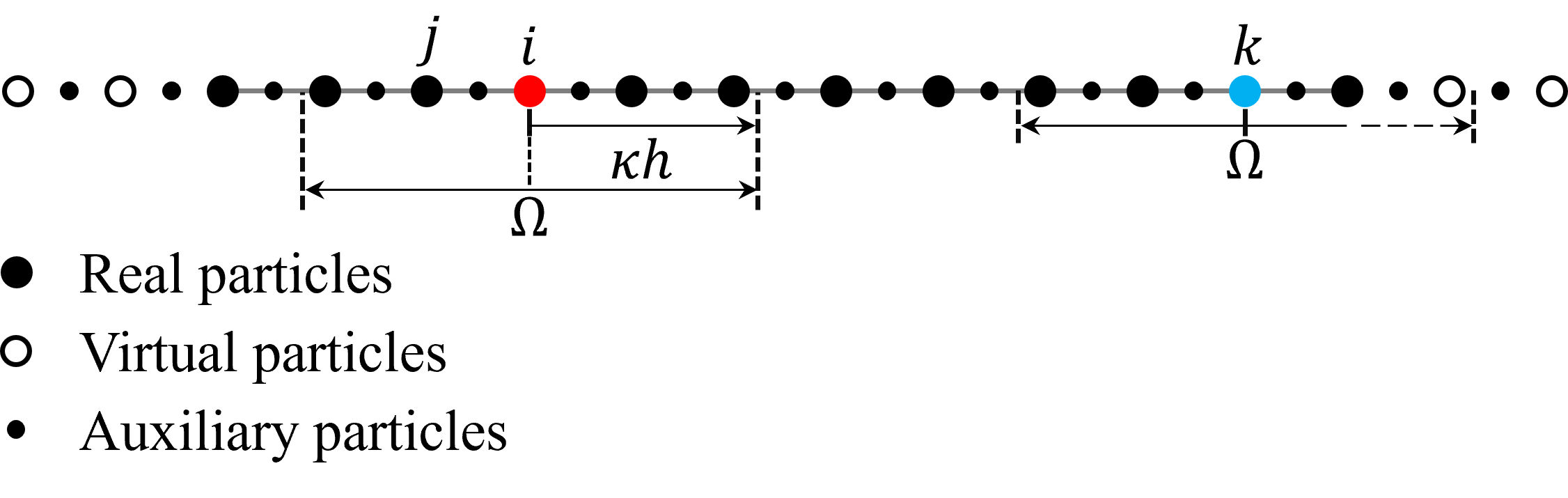}
    \caption{Schematic illustration of the discretized physical domain (containing real particles) and of the virtual domain (containing virtual particles). $\Omega$ is the support domain with radius $\kappa h$. Particle $i$ represents a particle with a fully supported domain, while particle $k$ represents a boundary-affected particle with a truncated support domain.}
    \label{fig:ِDomain}
\end{figure}

A second level of approximation occurs during the discretization process from continuum to particle representation. A brief summary of the discretization process is provided in the following \cite{ monaghan1988introduction,monaghan2005smoothed,liu2003smoothed}. Consider a generic continuous 1D domain that, in first approximation, can be thought of as a 1D beam. In SPH form, this domain can be discretized in a finite number of particles (Fig.~\ref{fig:ِDomain}) such that the generic particle $j$ corresponding to the finite domain of volume $V_j$ (in the initial continuous domain) around the particle has mass $m_j$ given by:

\begin{equation}
    m_j = V_j \rho_j \quad \text{for } j = 1, \dots, N
\end{equation}

where $N$ is the number of particles in which the continuous domain has been discretized. Based on this discretization approach, the infinitesimal volume $d\boldsymbol{x'}$ in the integral representation is replaced by the finite volume $V_j$ which leads to the following approximation of a generic field $f(\boldsymbol{x})$ (see Eq.~\ref{eq:SPH_int_approx}) at the particle located at $\boldsymbol{x}_i$:

\begin{equation}
    \langle f(\boldsymbol{x}_i) \rangle = \sum_{j=1}^N \frac{m_j}{\rho_j}f(\boldsymbol{x}_j)W_{ij}
    \label{eq:particle_discr}
\end{equation}

where $W_{ij}=W(\boldsymbol{x}_i-\boldsymbol{x}_j,h)$. Note that the discretization in Eq.~\ref{eq:particle_discr} does not limit the calculation to be performed only on the particle location $\boldsymbol{x}_i$, but it allows estimating the value of the function at any generic point $\boldsymbol{x}$ within the domain (even if not coincident with a particle); in this latter case the kernel function would take the more general form $W(\boldsymbol{x}-\boldsymbol{x}_j,h)$.

Following the above particle approximation, it is also immediate to obtain the discretized form of both the gradient and the integral of the function $f(\boldsymbol{x})$. The gradient of $f(\boldsymbol{x})$ will take the form:

\begin{equation}
     \langle \nabla f(\boldsymbol{x}_i) \rangle = \sum_{j=1}^N \frac{m_j}{\rho_j}f(\boldsymbol{x}_j) \nabla_i W_{ij}
     \label{eq:grad_xi}
\end{equation}

where $\nabla_i W_{ij}$ is the gradient of the kernel function $W_{ij}$ with respect to particle $i$.
Similarly, the discretized form of the integral (Eq.~\ref{eq:SPH_integral_kernel_final}) is given by:

\begin{equation}
     \langle \mathcal{I} (f(\boldsymbol{x}_i)) \rangle = \sum_{j=1}^N \frac{m_j}{\rho_j}f(\boldsymbol{x}_j) \Tilde{W}_{ij}
     \label{eq:Int_xi}
\end{equation}

where $\Tilde{W}_{ij}$ is the integral of the kernel function $W_{ij}$ that can be determined analytically. 
%The plot of this function for the cubic spline kernel is provided in Fig.~\ref{fig:Kernel}.

One of the complexities when using the SPH approximation is related to the particles whose support domain is truncated by the boundary (particle $k$ in Fig.~\ref{fig:ِDomain}). The corresponding reduction of the number of particles in the support domain leads to inconsistency in the approximated field at those particle locations. To address this issue, we introduce a virtual domain (discretized by means of virtual particles in Fig.~\ref{fig:ِDomain}) outside the real domain (discretized by real particles in Fig.~\ref{fig:ِDomain}) to ensure that there are enough particles in the support domain of the boundary-affected particles.  Moreover, to enhance the accuracy of approximation for all real particles, we also implement the kernel gradient correction. The corrected form of the kernel gradient proposed by Bonet and Lok \cite{bonet1999variational} is defined as:

\begin{equation}
     \nabla_{i}^{c} W_{ij} = \bm{\mathbb{L}}_{i}^{-1}\nabla_{i} W_{ij}
     \label{eq:corrected kernel gradient}
\end{equation}

with

\begin{equation}
     \bm{\mathbb{L}}_{i}^{-1} = -\sum_{j=1}^N \frac{m_j}{\rho_j} \nabla_i W_{ij}\otimes\boldsymbol{x}_{ij}
     \label{eq:L_norm correction}
\end{equation}

where $\boldsymbol{x}_{ij} = \boldsymbol{x}_i - \boldsymbol{x}_j$, and $\otimes$ indicates the tensor product. The kernel gradient correction satisfies first-order consistency, whereas the regular kernel gradient has zeroth-order consistency. Higher-order kernel gradient corrections have also been proposed in the literature \cite{asai2023class,nasar2021high}, but they generally demand higher computational cost. In this study, we employ the corrected form given in Eq. \ref{eq:corrected kernel gradient}.

%\subsection{A Modified SPH-based integration method}
\subsection{Improving the accuracy of the SPH-based integration scheme}\label{sec:SPH integral-modified}

This section introduces a modified form of the method designed to improve numerical performance and counteract some accuracy issues arising from the original SPH scheme. The utility of this modification will become evident in the following when numerical results will be presented.

The proposed improvement relies on a set of auxiliary particles introduced in the real domain and used to refine the SPH approximation. The use of auxiliary particles in SPH is not new and has been previously explored to enhance stability \cite{randles2000normalized,dyka1995approach,blanc2012stabilized}. As shown in Fig.~\ref{fig:ِDomain}, these auxiliary particles are placed at the midpoint between consecutive real particles, resembling the concept of midpoint selection commonly used in numerical integration. By defining these auxiliary particles, the SPH approximation of an integral at a generic particle $i$ is calculated using these auxiliary particles rather than the real particles in the domain. In other words, for a 1D domain, the SPH approximation of an integral given in Eq.~\ref{eq:Int_xi} is modified as:

\begin{equation}
     \langle \mathcal{I} (f({x}_i)) \rangle = \sum_{j=1}^N \frac{\bar{m}_j}{\bar{\rho}_j}f({\bar{x}}_j) \bar{\Tilde{W}}_{ij}
     \label{eq:Int_xi_auxiliary}
\end{equation}

where $\bar{x}_j = \frac{x_j+x_{j+1}}{2}$ is the location of an auxiliary particle between two adjacent real particles $j$ and $j+1$ and $\bar{\Tilde{W}}_{ij} = \Tilde{W}(x_i - \bar{x}_j,h)$ is the kernel integral which is calculated based on the distance between particle $i$ and an auxiliary particle 
$j$. Note that in SPH, the whole domain is represented by the real particles so the total mass and density are only allocated to these particles. Hence, to obtain the mass and density of the auxiliary particles (displayed as $\bar{m}_j$ and $\bar{\rho}_j$, respectively, in Eq.~\ref{eq:Int_xi_auxiliary}), we use the average value of the mass and density of the two adjacent real particles. That is to say, $\bar{m}_j = \frac{m_j+m_{j+1}}{2}$ and $\bar{\rho}_j = \frac{\rho_j+\rho_{j+1}}{2}$. 

Another important point concerns the evaluation of the integrand at an auxiliary particle, that is $f(\bar{x}_j)$. Although the focus of this study is on analytical functions that can be easily evaluated analytically, in a more practical application scenario the integrand may represent a physical field variable whose values are usually determined at the location of the real particles. In this study, to consider a more general approximation, we evaluate the function at an auxiliary particle as the average value of the function at two adjacent real particles, as shown in the following equation. 

\begin{equation}
     f(\bar{x}_j)=\frac{f(x_j)+f(x_{j+1})}{2}
     \label{eq:f-at-auxiliary}
\end{equation}

\subsection{Kernel function}

\begin{figure}[h!]
    \centering
    \includegraphics[width=0.5\textwidth, keepaspectratio]{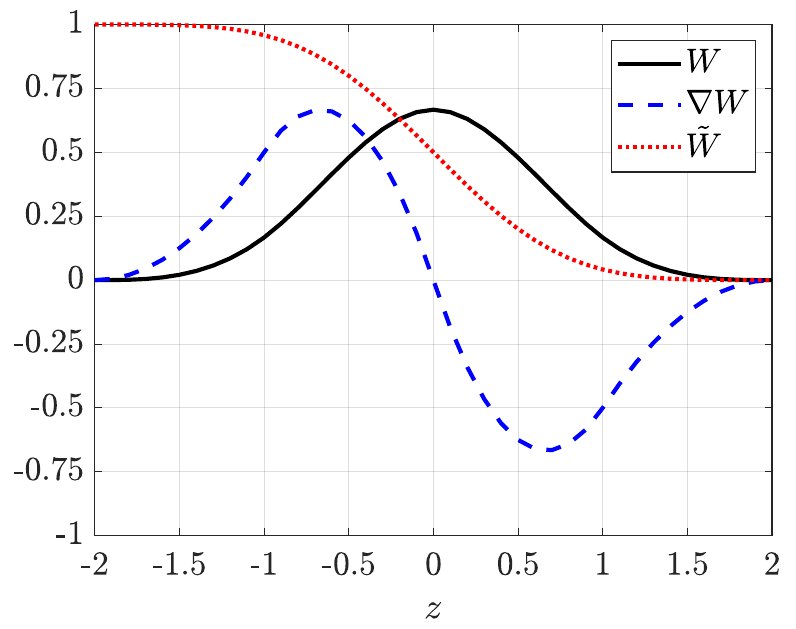}
    \caption{Plot of the cubic kernel function ($W$), its gradient ($\nabla W$), and its integral ($\Tilde{W}$).}
    \label{fig:Kernel}
\end{figure}

SPH can employ different types of kernel functions, each one displaying different properties \cite{liu2003smoothed}. In this study, we select the cubic spline function which is one of the most common kernels used in SPH approximation. Note that this specific choice of the kernel function will not affect the generality of the results. The cubic spline is a piecewise function defined as follows \cite{liu2003smoothed}:

\begin{equation}
    W(z,h)=C_d \cdot 
\begin{cases} 
      \frac{2}{3}-z^2+\frac{1}{2}z^3 & 0\leq z < 1 \\
      \frac{1}{6} (2-z)^3 & 1\leq z < 2 \\
      0 & z \geq 2 
   \end{cases}
\end{equation}

where $C_d$ is a constant whose value depends on the dimensionality of the problem. For a 1D problem, like the one considered in this study, $C_d=1/h$. Also, $z=\frac{r}{h}$ where $r$ is the absolute distance between two points at $\boldsymbol{x}$ and $\boldsymbol{x'}$. As evident from this equation, the radius of support domain $\Omega$ of this kernel function is $2h$ ($\kappa=2$ in Fig.~\ref{fig:ِDomain}). A plot of this kernel function, of its gradient, and its integral is provided in Fig.~\ref{fig:Kernel}.

As noted previously, to calculate the derivative and the integral of a generic function, the derivative and the integral of the kernel functions are needed. Due to the availability of the kernel functions in analytical form, both the derivative and the integral can be easily calculated analytically.

It should be mentioned that the generic point $\boldsymbol{x}$ in the SPH approximation represents the spatial location of a particle, which can include either one, two, or three coordinates depending on the dimensionality of the problem (i.e. 1D, 2D, or 3D). Since the focus of the present study is on a 1D problem, the regular (non-bold) symbol $x$ will be used in the rest of the paper.

\section{SPH approximation of fractional-order operators}\label{sec:SPH to FC}
Having introduced the main definitions for both fractional operators and SPH approximations, the calculation of fractional operators via the SPH approach can now be addressed. In general, it will be shown that the numerical calculation of fractional operators via the SPH method simply involves performing either integration or differentiation operations by exploiting the kernel functions, as discussed above.
%For the sake of brevity, only the SPH approximation of the left-handed fractional-order definitions are presented in the following subsections. The right-handed and Riesz-type counterparts can be obtained in a similar manner.

\subsection{CO operators} \label{sec:CO_SPH}
Starting from the definition of the RL fractional integral in Eq.~\ref{eq:RL_Integral_Left} and using the SPH representation of the integrated function described above, we arrive at the approximate SPH form of the RL integral:

\begin{equation}
    \langle \prescript{RL}{a}{\mathrm{I}}_{x_i}^{\alpha}f(x_i) \rangle = \frac{1}{\Gamma(\alpha)} \sum_{j=1}^N \frac{m_j}{\rho_j} \hat{f}(x_i,x_j) \Tilde{W}_{ij}\\
    \label{eq:SPH_RL_Integral_left}
\end{equation}

where
\begin{equation}
    \hat{f}(x_i,x_j) = (x_i-x_j)^{\alpha-1} f(x_j)
    \label{eq:f_hat}
\end{equation}

With the basic building block provided by the RL integral, obtaining the RL and Caputo derivatives simply requires performing an additional differentiation. For the RL derivative, following its definition in Eq.~\ref{eq:RL_Derivative_Left}, we perform a differentiation of order $n$ of the RL integral $\prescript{RL}{a}{\mathrm{I}}_x^{n-\alpha}$. The SPH approximation of $\prescript{RL}{a}{\mathrm{I}}_x^{n-\alpha}$  can be easily obtained from Eq.~\ref{eq:SPH_RL_Integral_left}. Applying further the SPH differentiation according to Eq.~\ref{eq:grad_xi}, we obtain the following expression:

\begin{equation} 
    \langle \prescript{RL}{a}{\mathrm{D}}_{x_i}^{\alpha}f(x_i) \rangle = \sum_{j=1}^N \frac{m_j}{\rho_j} \langle \prescript{RL}{a}{\mathrm{I}}_{x_i}^{n-\alpha}f(x_i) \rangle  \frac{d^nW_{ij}}{dx^n}  
    \label{eq:SPH_RL_derivative}
\end{equation}

The Caputo derivative can be calculated following a similar approach by performing, according to its definition in Eq.~\ref{eq:Caputo_CO_Derivative_Left}, first the derivative of the function and then taking the fractional integral:

\begin{equation}
    \langle \prescript{C}{a}{\mathrm{D}}_{x_i}^{\alpha}f(x_i) \rangle =  \frac{1}{\Gamma(n-\alpha)} \sum_{j=1}^N \frac{m_j}{\rho_j}\Tilde{f}(x_i,x_j) \Tilde{W}_{ij}
    \label{eq:SPH-Caputo_Derivative}
\end{equation}

where
\begin{equation}
    \Tilde{f}(x_i,x_j) = (x_i-x_j)^{n-\alpha-1} \langle \frac{d^n f(x^\prime)}{dx^n} \rangle
    \label{eq:f_tilde}
\end{equation}

and $\langle \frac{d^n f(x^\prime)}{dx^n} \rangle$ is the SPH approximation of the \textit{n}\textsuperscript{th}
 order derivative. In the present study, we consider $\alpha \in (0,1)$ and $ n = 1 $. So, $\frac{d^n}{dx^n}(\cdot)$ is simply $\frac{d}{dx}(\cdot)=\nabla(\cdot)$ which can be calculated using Eq.~\ref{eq:grad_xi} with the kernel gradient correction given in Eq.~\ref{eq:corrected kernel gradient}.

\subsection{VO operators}
For the case of VO operators, the approach follows an equivalent route to that outlined above. In Type-I VO considered in this study, where  $\alpha = \alpha (x)$, $\hat{f}(x_i,x_j)$ and $\tilde{f}(x_i,x_j)$ in Eqs.~\ref{eq:f_hat} and~\ref{eq:f_tilde} are redefined as:

\begin{equation}
    \hat{f}_{\scriptscriptstyle VO}(x_i,x_j) = \frac{1}{\Gamma\big(\alpha(x_i)\big)}(x_i-x_j)^{\alpha(x_i)-1} f(x_j)
    \label{eq:f_hat_VO}
\end{equation}

\begin{equation}
    \Tilde{f}_{\scriptscriptstyle VO}(x_i,x_j) = \frac{1}{\Gamma\big(n-\alpha(x_i)\big)}(x_i-x_j)^{n-\alpha(x_i)-1} \langle \frac{d^n f(x^\prime)}{dx^n} \rangle
    \label{eq:f_tilde_VO}
\end{equation}

which lead to the SPH approximations of the VO RL integral, the VO RL derivative, and the VO Caputo derivative, respectively, as follows:

\begin{equation}
    \langle \prescript{RL}{a}{\mathrm{I}}_{x_i}^{\alpha(x_i)}f(x_i) \rangle \simeq \sum_{j=1}^N \frac{m_j}{\rho_j} \hat{f}_{\scriptscriptstyle VO}(x_i,x_j) \Tilde{W}_{ij} 
    \label{eq:SPH_VO_RL_Integral_left}   
\end{equation}

\begin{equation}
    \langle \prescript{RL}{a}{\mathrm{D}}_{x_i}^{\alpha(x_i)}f(x_i) \rangle  \simeq \sum_{j=1}^N \frac{m_j}{\rho_j} \langle \prescript{RL}{a}{\mathrm{I}}_{x_i}^{n-\alpha(x_i)}f(x_i) \rangle  \frac{d^nW_{ij}}{dx^n} 
    \label{eq:SPH_VO_RL_derivative_left}   
\end{equation}

\begin{equation}
    \langle \prescript{C}{a}{\mathrm{D}}_{x_i}^{n-\alpha(x_i)}f(x_i) \rangle \simeq \sum_{j=1}^N \frac{m_j}{\rho_j} \Tilde{f}_{\scriptscriptstyle VO}(x_i,x_j) \Tilde{W}_{ij} 
    \label{eq:SPH_VO_Caputo_derivative_left}   
\end{equation}

The SPH integral approximations presented in this section are based on the standard form presented in Eq.~\ref{eq:Int_xi}. To obtain an approximation using modified SPH integration given in Eq.~\ref{eq:Int_xi_auxiliary}, one simply needs to use Eq.~\ref{eq:Int_xi_auxiliary} in place of Eq.~\ref{eq:Int_xi}.

\section{Addressing the kernel singularity}\label{sec:singularity}
The numerical solution of fractional-order operators is well-known to pose challenges due to the singularity of the power law kernel at the domain boundary. This topic has been the object of many studies, and a variety of approaches have been proposed over the years. Examples include, but are not limited to, 
an analytical approach based on integration by parts \cite{jacobs2023order}, a semi-analytical approach based on the Hadamard's finite-part integral\cite{delbourgo1994approximate,elliott1995three,kaya1987solution}, and numerical integration schemes specifically conceived for singular systems of equations  \cite{korsunsky1998gauss,ioakimidis1989location,carley2007numerical}.

In this study, we choose an analytical approach that eliminates the singularity by performing a common manipulation of the definition of the fractional-order operator by means of integration by parts \cite{jacobs2023order}. Starting from the definition of the left-handed CO RL integral in Eq.~\ref{eq:RL_Integral_Left}, and considering $0<\alpha<1$, the application of integration by parts yields:

\begin{equation}
     \prescript{RL}{a}{\mathrm{I}}_x^{\alpha}f(x) = \frac{1}{\Gamma(\alpha)} \left(\frac{f(a) (x-a)^{\alpha}}{
     \alpha} +\int_a^x f'(x') \frac{(x-x')^{\alpha}}{\alpha} dx' \right)
    \label{eq:RL_Integral_nonsing_left}      
\end{equation}

It is evident that the kernel singularity is removed and the expression can be readily integrated. At the same time, from Eq.~\ref{eq:RL_Integral_nonsing_left} it is seen that the price to pay to eliminate the singularity is to calculate the derivative $f'(x')$. This means that, unless the function is known and the derivative can be calculated analytically, the approach requires numerical differentiation that ultimately increases the computational cost. In the case of SPH, the estimate of a derivative is reduced to performing a simple algebraic operation where the analytically differentiated kernel is multiplied by the function of interest. However, the operation still requires a summation over the particles in the domain of interest, which does increase the computational cost of the algorithm. 

While this approach will be used in the numerical results section, we mention that other strategies could also be explored to avoid the burden due to the evaluation of the additional derivative. The non-singular form presented in Eq.~\ref{eq:RL_Integral_nonsing_left} is also valid for Type-I VO. In this case, since $\alpha =\alpha(x)$, it remains independent of the integration variable $x'$ and can be treated as a constant within the integrand of the fractional integral definition. However, for the same reason, this equation does not apply to Type-II and Type-III VO, where $\alpha$ depends on the variable $x'$.  

From the non-singular form of the RL integral in Eq.~\ref{eq:RL_Integral_nonsing_left}, it is immediate to obtain the non-singular form of the RL derivative:
%\begin{comment}
\begin{equation}
\begin{split}
     \prescript{RL}{a}{\mathrm{D}}_x^{\alpha}f(x) &=\frac{d}{dx} \left( \prescript{RL}{a}{\mathrm{I}}_x^{1-\alpha}f(x) \right) \\&=\frac{d}{dx} \left[ \frac{1}{\Gamma(1-\alpha)} \left(\frac{f(a) (x-a)^{1-\alpha}}{
     1-\alpha} +\int_a^x f'(x') \frac{(x-x')^{1-\alpha}}{1-\alpha} dx' \right)  \right] 
    \label{eq:RL_Derivative_nonsing}      
\end{split}
\end{equation}
%\end{comment}

This equation applies to both CO and Type-I VO. In the case of CO, it simplifies further and the final form, also presented in \cite{diethelm2010analysis}, is given by:

\begin{equation}
     \prescript{RL}{a}{\mathrm{D}}_x^{\alpha}f(x) = \frac{1}{\Gamma(1-\alpha)} \left(\frac{f(a)}{(x-a)^\alpha} +\frac{d}{dx}\int_a^x f'(x') \frac{(x-x')^{1-\alpha}}{1-\alpha} dx' \right)
    \label{eq:RL_Derivative_nonsing_CO}      
\end{equation}

As seen in Eq.~\ref{eq:RL_Derivative_nonsing_CO}, although the singularity is eliminated at $x = x'$, it still remains at the boundary $x = a$ when $f(a) \neq 0$.

%\KG{in this case, there is singularity at $x = a$ (from the first term). Just when $f(a) = 0$, we have $0/0$ where in the limit results in $0$. For example, for $f(x) = sin(x)$ and $a = 0$, we have the limit case $\lim_{x \to 0} \frac{0}{x^\alpha } = 0 $ ($\alpha>0$)}

In a similar way, the non-singular form of the Caputo derivative can be  obtained by integrating by parts Eq.~\ref{eq:Caputo_CO_Derivative_Left}. As also reported in \cite{jacobs2023order}, the resulting form is:

\begin{equation}
\begin{split}
     \prescript{C}{a}{\mathrm{D}}_x^{\alpha}f(x) &= \frac{1}{\Gamma(1-\alpha)} \left(\frac{f'(a) (x-a)^{1-\alpha}}{1-\alpha} + \int_a^x f''(x') \frac{(x-x')^{1-\alpha}}{1-\alpha} dx' \right) \\
     %&= \frac{1}{\Gamma(2-\alpha)} \left(\frac{f'(a) }{(x-a)^{\alpha-1}} + \int_a^x  \frac{f''(x')}{(x-x')^{\alpha-1}} dx' \right)
    \label{eq:Caputo_Derivative_nonsing_left}      
\end{split}
\end{equation}

As previously discussed for the non-singular form of RL integral (Eq.~\ref{eq:RL_Integral_nonsing_left}), the same reasoning ensures that this equation is also valid for both CO and Type-I VO. Moreover, as already observed for the RL operators, the order of the derivative appearing in the operator increases, but the singularity is eliminated.

The SPH approximation of the non-singular expressions can be derived by following the same procedure outlined for the original fractional-order operators. As a complementary equation to the SPH approximations,  
the second-order derivative required to define the non-singular form of the Caputo derivative (see Eq.~\ref{eq:Caputo_Derivative_nonsing_left}) can be approximated using a formula proposed by Brookshaw \cite{brookshaw1985method}, as follows:

\begin{equation}
     \langle \nabla^2 f(x_i) \rangle \simeq -2\sum_{j=1}^N \frac{m_j}{\rho_j}\big(f(x_j)-f(x_i)\big)\frac{\boldsymbol{x}_{ij}.\nabla_iW_{ij}}{|\boldsymbol{x}_{ij}|^2+\eta^2}
    \label{eq:Brookshaw}      
\end{equation}

 where $\eta=0.0001(1/h)^2$ is a parameter introduced to prevent a singularity. Note that in this alternative expression, the second-order derivative is approximated using the gradient of the kernel function (rather than its Laplacian) whose corrected form was defined in Eq.~\ref{eq:corrected kernel gradient}. 

\section{Numerical results} \label{sec: results}
This section presents the results of a numerical analysis to investigate the performance of the SPH approximation of fractional-order operators. %Depending on the form of the operator, results are validated via either analytical or numerical results. 
The SPH implementation of the fractional-order operators presented above is general and can be applied to virtually any function. In the following, we choose specific examples based on the availability of their analytical solutions so to perform a direct and independent validation. Four basic functions that have been selected are the trigonometric functions $\sin(\pi x)$ and $\cos(\pi x)$, the exponential function $\exp(x)$, and the shifted polynomial function $(x-1)^3$. The analytical solutions of these functions for the (left-handed) Caputo derivative, RL derivative, and the RL integral are presented in \ref{sec:Appendix A}. These solutions are either available in the literature \cite{shchedrin2018exact,herrmann2011fractional} or can be obtained based on hypergeometric functions via commercial software such as Wolfram Mathematica \cite{Mathematica}.

To assess the accuracy of the results produced via SPH approximation, the absolute error, the relative \textit{$L_2$}-error, and the \textit{$R^2$}-score will be reported for each numerical evaluation. The absolute error as a point-wise metric error is defined as: 
%the absolute difference between the analytical and the approximated solutions:

\begin{equation}
    \text{Abs. Error} = |f_i-\langle f \rangle_i|
    \label{eq:Abs Error}
\end{equation}

where $f_i$ and $\langle f \rangle_i$ are the analytical and the approximate values calculated at particle $i$. The relative \textit{$L_2$}-error quantifies the average relative error and is expressed as \cite{nair2024multiple}:

\begin{equation}
    L_2 = \frac{\left(\sum_{i=1}^{N} |f_i-\langle f \rangle_i|^2\right)^{1/2}}{\left(\sum_{i=1}^{N} |f_i|^2\right)^{1/2}}
    \label{eq:L2 norm}
\end{equation}

 where $N$ indicates the total number of particles. A lower relative \textit{$L_2$}-error corresponds to higher accuracy in the approximated solutions.

The \textit{$R^2$}-score, also known as the coefficient of determination, measures how well the approximated values match the analytical values, as given by the following equation \cite{nair2024multiple}:

\begin{equation}
    R^2 = 1 - \frac{\sum_{i=1}^{N} (f_i-\langle f \rangle_i)^2}{\sum_{i=1}^{N} (f_i - \bar{f})^2}
    \label{eq:R2 score}
\end{equation}

where $\bar{f}$ represents the average of the analytical solutions. The maximum possible $R^2$-score value is $1$ which occurs when the approximate and analytical solutions coincide. Therefore, the higher the $R^2$ value the higher the accuracy.

Concerning the results presented in the following, all SPH approximations were calculated using a particle spacing $s=0.0125$, with a total of $N=401$ particles and a smoothing length of $h=1.1s$. The fractional-order was set at $\alpha=0.75$ for CO and $\alpha(x) = 0.5 + 0.3\sin(4\pi x)$ for Type-I VO that varies between $0$ and $1$. All calculations were performed via the commercial software MATLAB R2024b.

\subsection{CO operators}
We first test the SPH approximation of the CO operators provided in \S\ref{sec:CO_SPH}. Overall, results will be presented in figures following a two-column arrangement. Plots in the left column (i.e. (a1), (b1), (c1), and (d1)) compare the SPH approximation with the analytical solution for the four selected functions. On the other hand, plots in the right column (i.e. (a2), (b2), (c2), and (d2)) present the corresponding point-wise absolute error along with the relative \textit{$L_2$}-error and $R^2$-score listed in the inset. 

\begin{figure}[h!]
    \centering
    \includegraphics[width=0.8\textwidth, keepaspectratio]{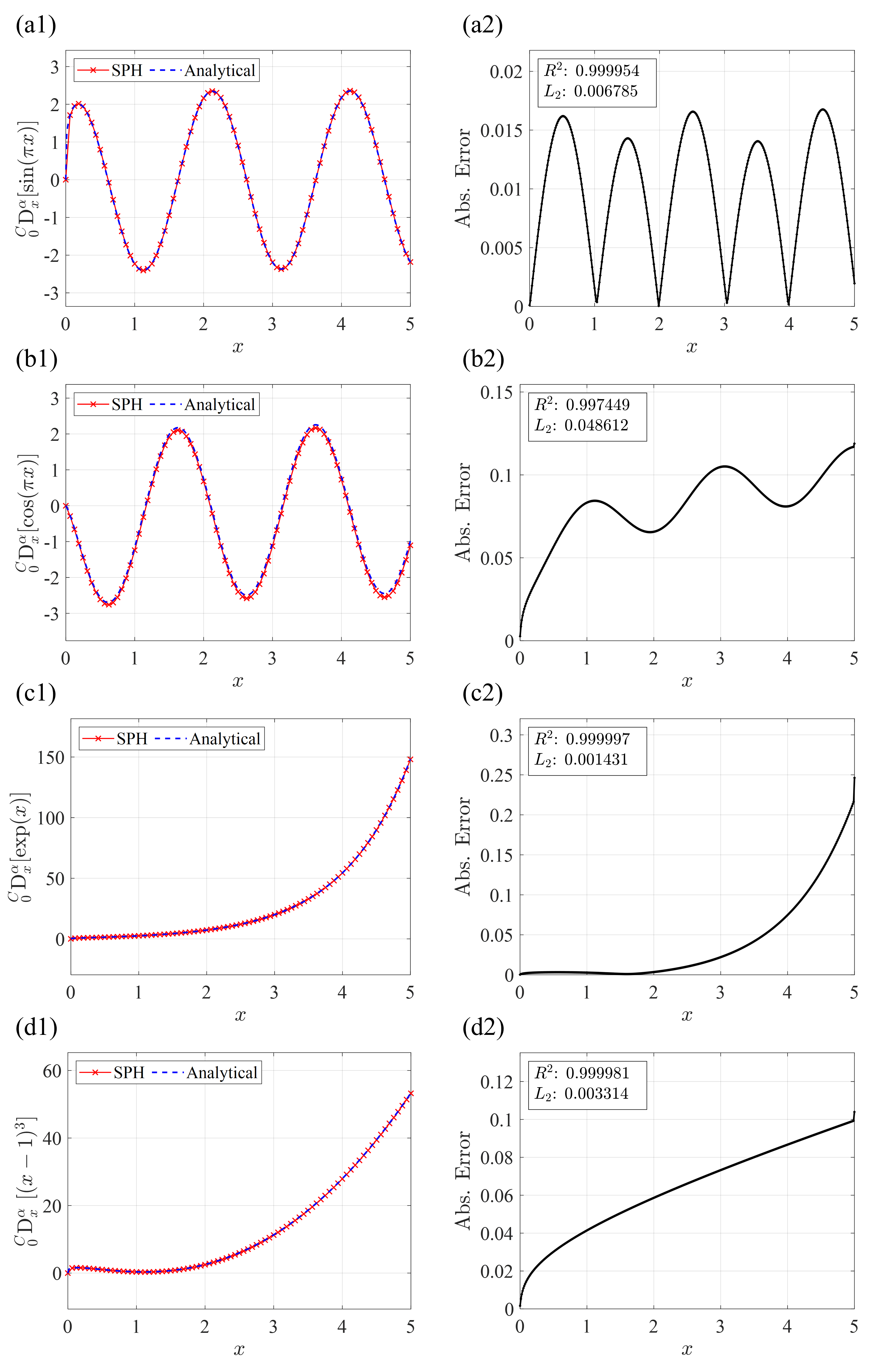}
    \caption{Left-handed CO Caputo derivative: Comparison between the SPH approximation and the analytical solution for (a1) $\sin (\pi x)$, (b1) $\cos (\pi x)$, (c1) $\exp(x)$, and (d1) $(x-1)^3$. The right column ((a2), (b2), (c2), and (d2)) shows the corresponding absolute error, including the $L_2$-error and the $R^2$-score metrics.}
    \label{fig:left-Caputo-CO}
\end{figure}

\begin{figure}[h!]
    \centering
    \includegraphics[width=0.8\textwidth, keepaspectratio]{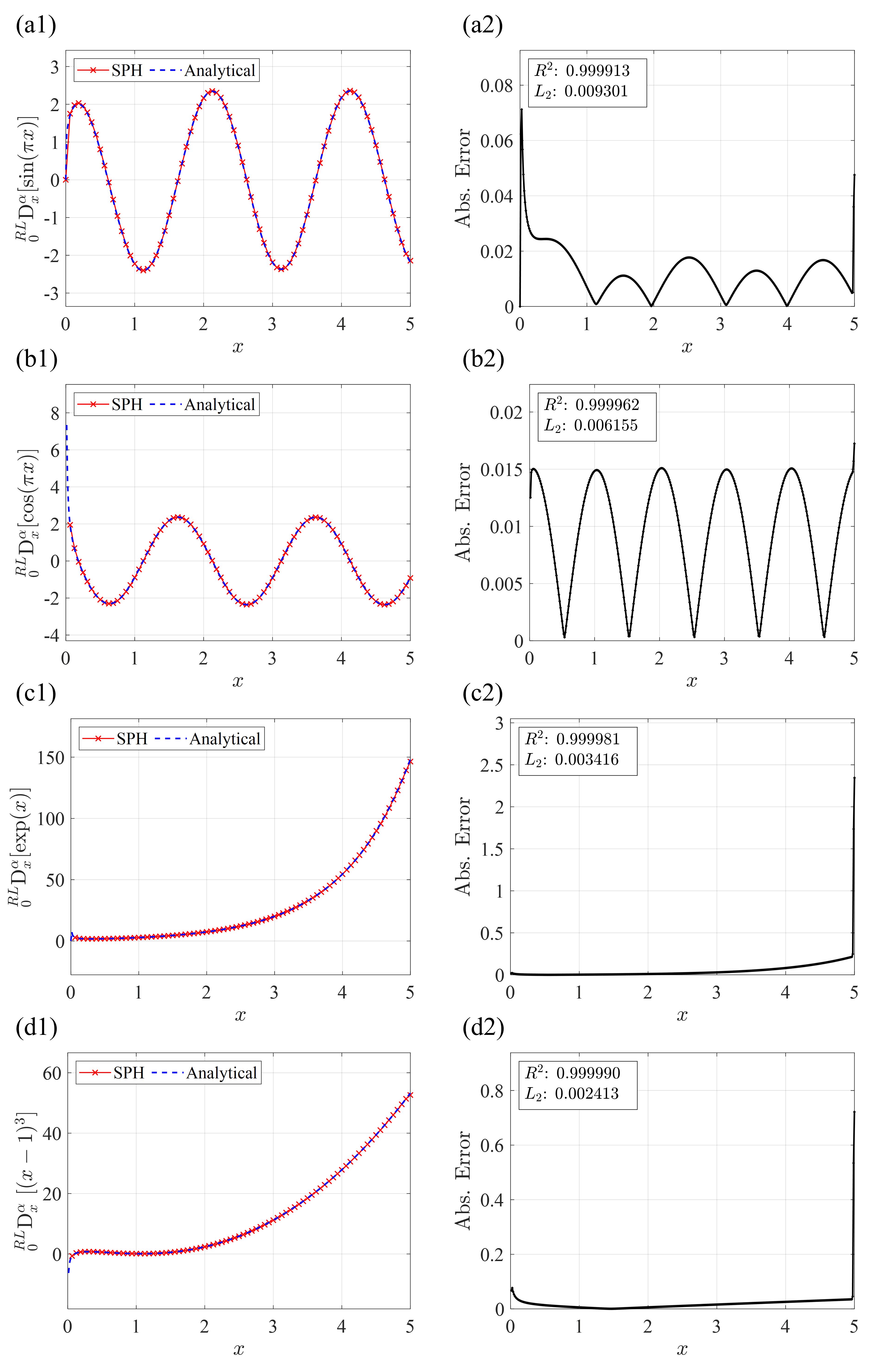}
    \caption{Left-handed CO RL derivative: Comparison between the SPH approximation and the analytical solution for (a1) $\sin (\pi x)$, (b1) $\cos (\pi x)$, (c1) $\exp(x)$, and (d1) $(x-1)^3$. The right column ((a2), (b2), (c2), and (d2)) shows the corresponding absolute error, including the $L_2$-error and the $R^2$-score metrics.}
    \label{fig:left-RL derivative-CO}
\end{figure}

\begin{figure}[h!]
    \centering
    \includegraphics[width=0.8\textwidth, keepaspectratio]{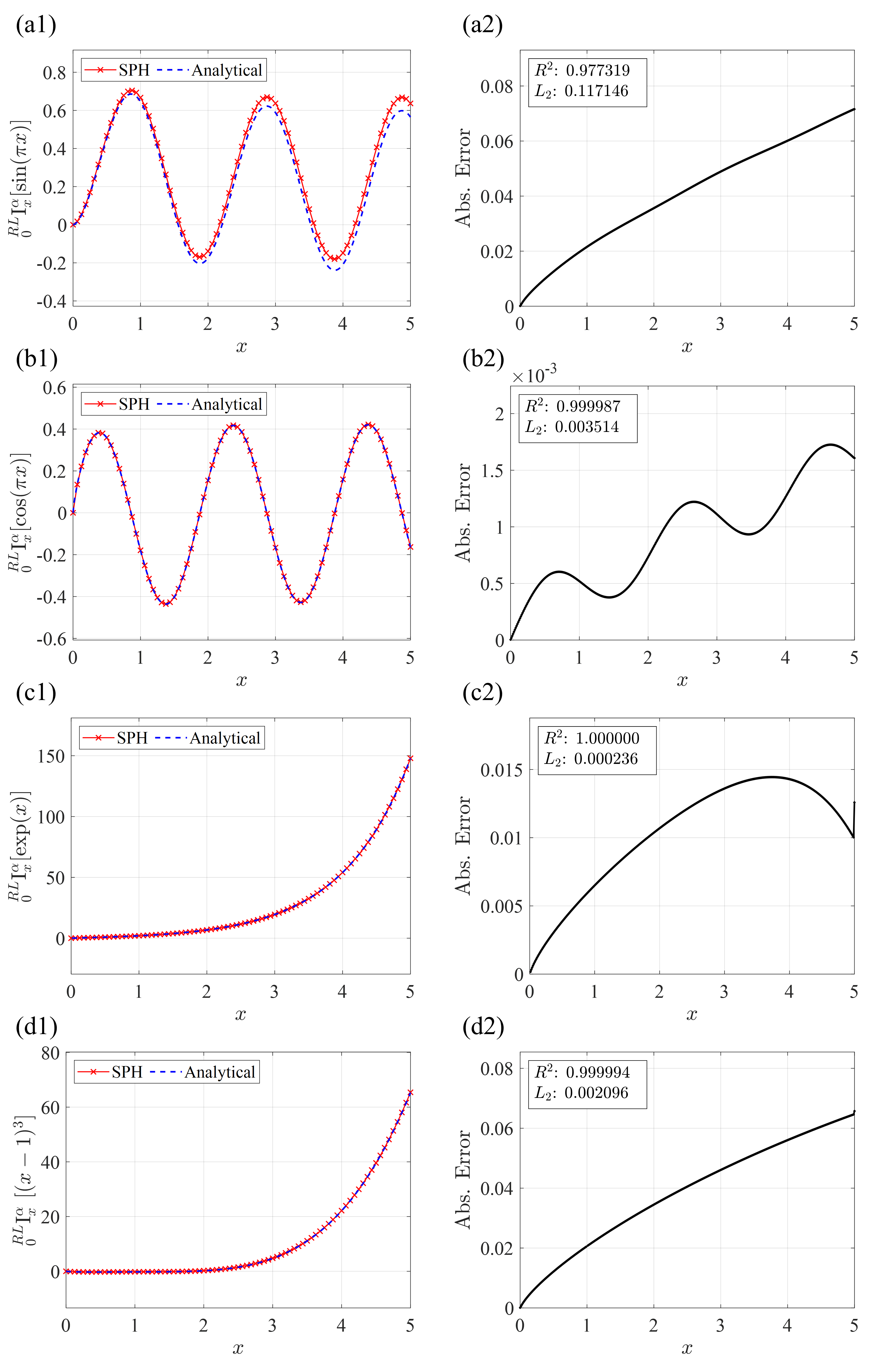}
    \caption{Left-handed CO RL integral: Comparison between the SPH approximation and the analytical solution for (a1) $\sin (\pi x)$, (b1) $\cos (\pi x)$, (c1) $\exp(x)$, and (d1) $(x-1)^3$. The right column ((a2), (b2), (c2), and (d2)) shows the corresponding absolute error, including the $L_2$-error and the $R^2$-score metrics.}
    \label{fig:left-RL integral-CO}
\end{figure}

\begin{figure}[h!]
    \centering
    \includegraphics[width=0.8\textwidth, keepaspectratio]{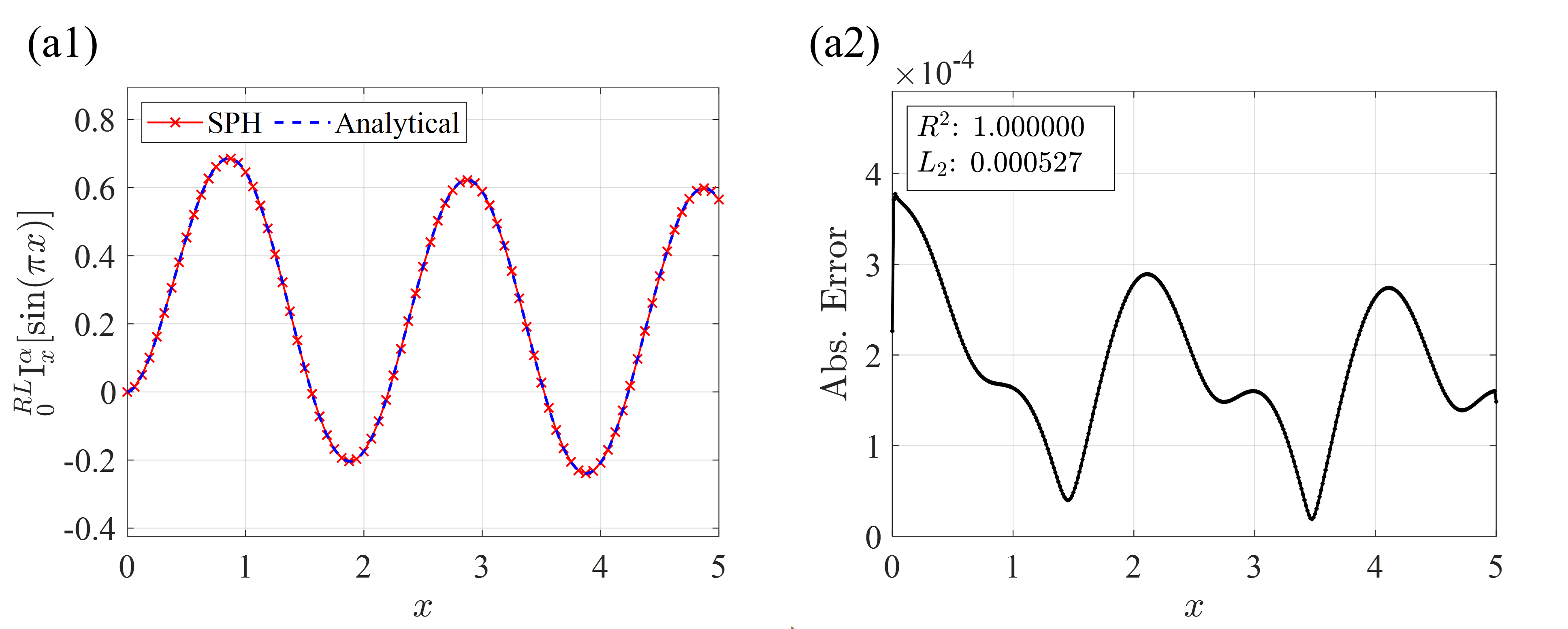}
    \caption{Left-handed CO RL integral: (a1) Comparison between the modified SPH approximation and the analytical solution for $\sin (\pi x)$. (a2) Corresponding absolute error, including the $L_2$-error and the $R^2$-score metrics.}
    \label{fig:left-RL integral-CO-modified}
\end{figure}

More specifically, Fig.~\ref{fig:left-Caputo-CO} displays the results of the Caputo derivative. A good agreement is observed between the approximate and the analytical solution, with the least accurate estimation corresponding to the function $ \cos (\pi x)$ (plots (b1) and (b2)) which displays $L_2=0.048612$ and $R^2=0.997449$. Fig.~\ref{fig:left-RL derivative-CO} presents the results for the RL derivative, highlighting the accuracy of the SPH approximations. The lowest accuracy is observed for the function $f(x) = \sin (\pi x)$ (plots (a1) and (a2)), with a maximum relative \textit{$L_2$}-error of $0.009301$ and the minimum $R^2$-score of $0.999913$. Similarly, the findings for RL integral are given in Fig.~\ref{fig:left-RL integral-CO}. Contrarily to the results obtained for both the Caputo and the RL derivatives, the error in the SPH approximation of $\prescript{RL}{0}{\mathrm{I}}_x^{\alpha}[\sin(\pi x)]$ (plots (a1) and (a2)) is relatively high. For the other three test functions, results are aligned with the analytical values. 
In an attempt to counteract this error and increase accuracy, 
we apply the modified form of the SPH integral approximation developed in \S\ref{sec:SPH integral-modified}, to the calculation of $\prescript{RL}{0}{\mathrm{I}}_x^{\alpha}[\sin(\pi x)]$. As illustrated in  Fig.~\ref{fig:left-RL integral-CO-modified}, the results obtained using this modified approach show a significant reduction in error and improved accuracy. More specifically, the \textit{$L_2$}-error decreases from $0.117146$ to $0.000527$ and the $R^2$-score increases from $0.977319$ to $1.000000$.

\subsection{VO operators}\label{sec:Results-VO}
\begin{figure}[h!]
    \centering
    \includegraphics[width=0.8\textwidth, keepaspectratio]{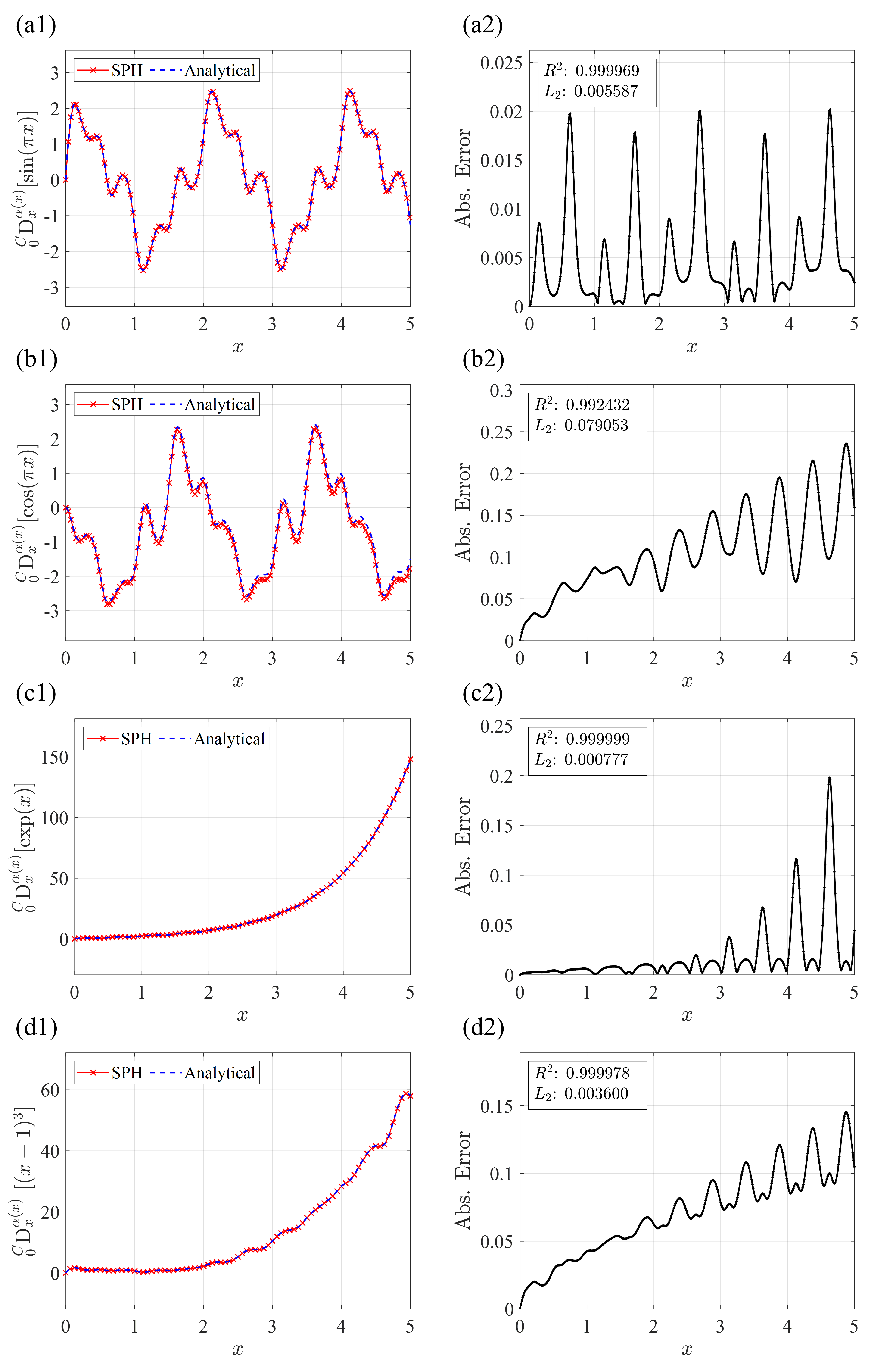}
    \caption{Left-handed Type-I VO Caputo derivative: Comparison between the SPH approximation and the analytical solution for (a1) $\sin (\pi x)$, (b1) $\cos (\pi x)$, (c1) $\exp(x)$, and (d1) $(x-1)^3$. The right column ((a2), (b2), (c2), and (d2)) shows the corresponding absolute error, including the $L_2$-error and the $R^2$-score metrics.}
    \label{fig:left-Caputo-VO}
\end{figure}

\begin{figure}[h!]
    \centering
    \includegraphics[width=0.8\textwidth, keepaspectratio]{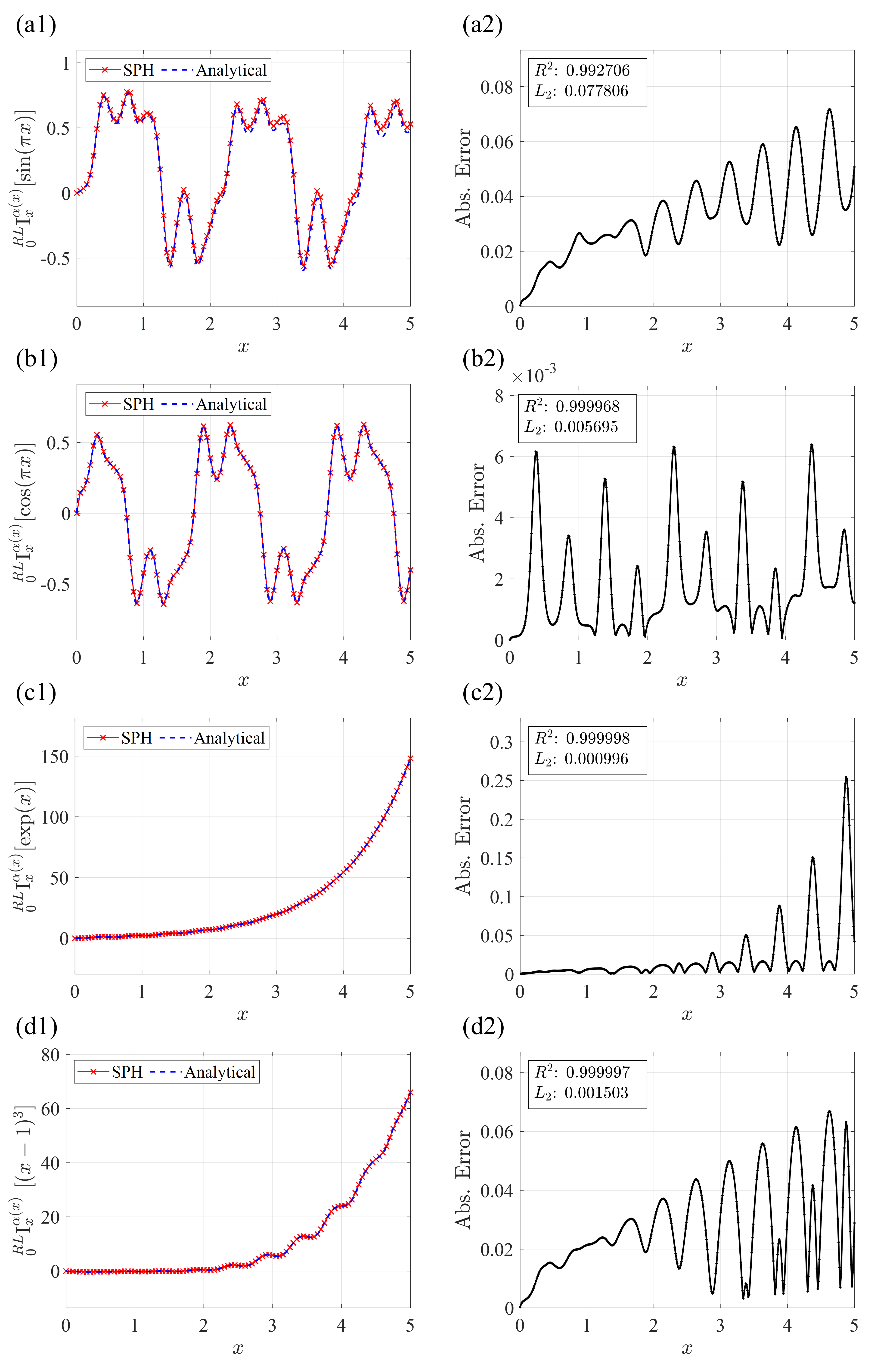}
    \caption{Left-handed Type-I VO RL integral: Comparison between the SPH approximation and the analytical solution for (a1) $\sin (\pi x)$, (b1) $\cos (\pi x)$, (c1) $\exp(x)$, and (d1) $(x-1)^3$. The right column ((a2), (b2), (c2), and (d2)) shows the corresponding absolute error, including the $L_2$-error and the $R^2$-score metrics.}
    \label{fig:left-RL integral-VO}
\end{figure}

\begin{figure}[h!]
    \centering
    \includegraphics[width=0.8\textwidth, keepaspectratio]{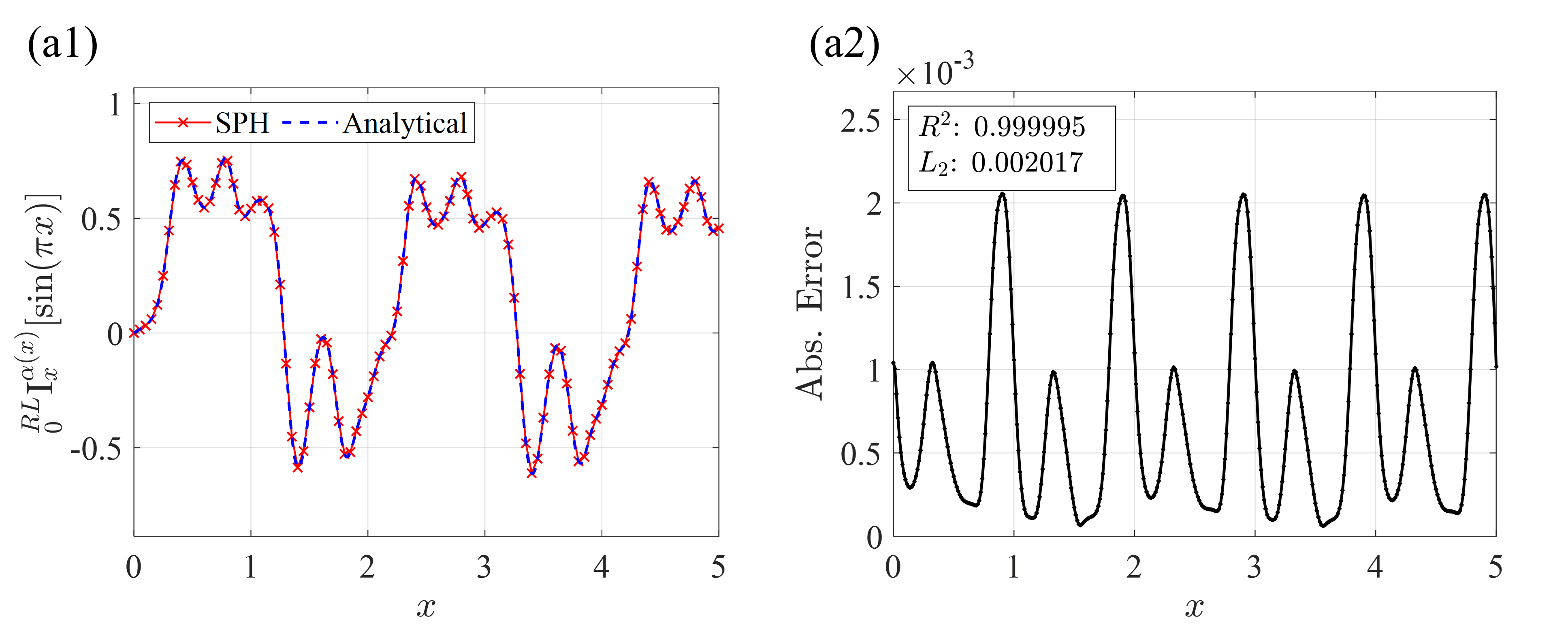}
    \caption{Left-handed VO RL integral: (a1) Comparison between the modified SPH approximation and the analytical solution for $\sin (\pi x)$. (a2) Corresponding absolute error, including the $L_2$-error and the $R^2$-score metrics.}
    \label{fig:left-RL integral-VO-modified}
\end{figure}

In a similar way, we investigated the performance of the proposed SPH approximation algorithms for the Type-I VO definition of the fractional-order operators. As discussed in \S\ref{sec:VO}, the available analytical solutions for the CO Caputo derivative and the RL integral also hold for their Type-I VO counterparts. Therefore, in the following, we present results for only these two types of VO operators. Fig.~\ref{fig:left-Caputo-VO} compares the SPH approximation with the analytical solution of the VO Caputo derivative. Across all chosen sample functions, the numerical results closely match the analytical solutions with minimal deviations. Similar to the case of the CO Caputo derivative, the worst-case performance of the SPH estimations is observed for the function $\cos(\pi x)$ with the maximum $L_2$-error of 0.079053 and minimum $R^2$-score of 0.992432 (plots (b1) and (b2)). 
The accuracy of the SPH approximation for the VO RL integral is demonstrated in Fig. \ref{fig:left-RL integral-VO}. The numerical results show strong agreement with the analytical solution for all functions. Notably, the largest deviation occurs for the function $\sin(\pi x)$ where the maximum $L_2$-error reaches 0.077806 and the minimum $R^2$-score is 0.992706 (plots (a1) and (a2)). This represents the lowest accuracy among the tested functions, yet the overall accuracy remains high, highlighting the reliability of the SPH method. To show the performance and applicability of the modified SPH integral approximation for VO operators, we applied this modified form to the calculation of $\prescript{RL}{0}{\mathrm{I}}_x^{\alpha(x)}[\sin(\pi x)]$. The results, presented in Fig~\ref{fig:left-RL integral-VO-modified}, confirm the effectiveness of this approach when using VO operators. 

\section*{Conclusions} \label{sec: conclusion}
In the present study, we developed a SPH-based formalism for the numerical evaluation of fractional-order operators. A few of the most common definitions of fractional operators were explicitly addressed although the formulation is general and applicable to any of the many definitions of fractional operators. %including left-handed, right-handed, and Riesz-type definitions of both constant-order (CO) and variable-order (VO) Riemann-Liouville (RL) integral, RL derivative, and  Caputo derivative. 
Given the differ-integral nature of these operators, our approach includes both SPH-based differentiation and integration. The calculation of integrals via the SPH method was performed according to two different methods. The first method followed the standard SPH definition that approximates the value of the operators at the locations of real particles. The second method involved the introduction of a set of auxiliary particles within the domain in order to increase the accuracy of the approximation.
To validate the theoretical formulation and the effectiveness of the resulting numerical method, we applied it to different types of functions including trigonometric, exponential, and polynomial functions. The direct comparison of the numerical results with the corresponding analytical solutions confirmed the viability and accuracy of the method. The proposed SPH-based formulation to numerically evaluate fractional operators lays the necessary foundation to solve fractional-order differential equations as well as general integral equations.\\

\noindent \textbf{Competing interests}

The authors declare no competing interest.\\

\noindent \textbf{Acknowledgments}
The authors gratefully acknowledge the partial financial support of the Army Research Office under the project  \# 82248-TE.

\renewcommand{\theequation}{A\arabic{equation}}
\setcounter{equation}{0}

\clearpage
\appendix
\renewcommand{\thesection}{Appendix}
\section{} \label{sec:Appendix A}
In this appendix, we provide analytical solutions for CO operators with $0<\alpha<1$, which were used to generate the comparison for the SPH numerical results. As discussed in \S\ref{sec:Results-VO}, the analytical solutions of CO Caputo derivative and CO RL integral are applicable to the Type-I VO of these fractional operators.

\subsection*{Left-handed Caputo derivative\cite{shchedrin2018exact,herrmann2011fractional}}

\vspace{10pt}
\begin{equation}
\begin{split}
    &\prescript{C}{0}{\mathrm{D}}_x^{\alpha}[\sin(\beta x)] = \frac{\beta x^{1 - \alpha}}{\Gamma(2-\alpha)} 
    {}_{1}F_{2} \left[ 1; 1-\frac{\alpha}{2},\frac{3}{2}- \frac{\alpha}{2}; -\frac{(\beta x)^{2}}{4} \right]
    \label{eq:Caputo Analytical sin}      
\end{split}
\end{equation}

\begin{equation}
\begin{split}
    &\prescript{C}{0}{\mathrm{D}}_x^{\alpha}[\cos(\beta x)] =- \frac{\beta^2 x^{2 - \alpha}}{\Gamma(3-\alpha)} 
    {}_{1}F_{2} \left[ 1; \frac{3}{2}-\frac{\alpha}{2},2- \frac{\alpha}{2}; -\frac{(\beta x)^{2}}{4} \right]
    \label{eq:Caputo Analytical cos}      
\end{split}
\end{equation}

% From Herrmann book
\begin{equation}
\begin{split}
    \prescript{C}{0}{\mathrm{D}}_x^{\alpha}[\exp(\beta x)] = \operatorname{sgn}(x)\left[(\beta\operatorname{sgn}(x)\right] ^\alpha\exp(\beta x)\Bigg( 1-\frac{\Gamma(1-\alpha,\beta x)}{\Gamma(1-\alpha)} \Bigg) 
    \label{eq:Caputo Analytical exp}      
\end{split}
\end{equation}

% From "Exact results for a fractional derivative of elementary functions" paper
\begin{equation}
\begin{split}
    &\prescript{C}{0}{\mathrm{D}}_x^{\alpha}[(x+\beta)^n] = \frac{\beta^{n-1}}{x^{\alpha-1}} \frac{n}{\Gamma(2-\alpha)}
    {}_{2}F_{1} \left[ 1,1-n;2-\alpha; -\frac{x}{\beta} \right]
    \label{eq:Caputo Analytical shifted poly}      
\end{split}
\end{equation}

where ${}_{p}F_{q} (a_1,\dots,a_p;b_1,\dots,b_p;z)$ is the generalized hypergeometric function and $\Gamma(\nu,\gamma)$ is the incomplete gamma function \cite{gradshteyn2014table}.
\subsection*{Left-handed RL derivative\cite{herrmann2011fractional,Mathematica}}
\vspace{10pt}
% from Mathematica built-in function
\begin{equation}
\begin{split}
    &\prescript{RL}{0}{\mathrm{D}}_x^{\alpha}[(\sin(\beta x)] = \frac{\beta x^{1 - \alpha}2^{\alpha - 1}\sqrt{\pi}}{\Gamma(1-\frac{\alpha}{2})\Gamma(\frac{3}{2}-\frac{\alpha}{2})}{}_{1}F_{2} \left[ 1; 1-\frac{\alpha}{2},\frac{3}{2}-\frac{\alpha}{2}; -\frac{(\beta x)^{2}}{4} \right]
    \label{eq:RL derivative Analytical sin}      
\end{split}
\end{equation}

% from Mathematica built-in function
\begin{equation}
\begin{split}
    &\prescript{RL}{0}{\mathrm{D}}_x^{\alpha}[(\cos(\beta x)] = \frac{ x^{- \alpha}2^{\alpha}\sqrt{\pi}}{\Gamma(\frac{1}{2}-\frac{\alpha}{2})\Gamma(1-\frac{\alpha}{2})}{}_{1}F_{2} \left[ 1; \frac{1}{2}-\frac{\alpha}{2},1-\frac{\alpha}{2}; -\frac{(\beta x)^{2}}{4} \right]
    \label{eq:RL derivative Analytical cos}      
\end{split}
\end{equation}

% From Herrmann book
\begin{equation}
\begin{split}
    \prescript{RL}{0}{\mathrm{D}}_x^{\alpha}[(\exp(\beta x)] = \operatorname{sgn}(x)\left[(\beta\operatorname{sgn}(x)\right] ^\alpha\exp(\beta x)\Bigg( 1-\frac{\Gamma(-\alpha,\beta x)}{\Gamma(-\alpha)} \Bigg) 
    \label{eq:RL derivative Analytical exp}      
\end{split}
\end{equation}

% from Mathematica built-in function
\begin{equation}
\begin{split}
    &\prescript{RL}{0}{\mathrm{D}}_x^{\alpha}[(x+\beta)^n] = \frac{(\beta+x)^{n}}{x^{\alpha}}\Big(\frac{\beta+x}{\beta}\Big)^{-n} \frac{1}{\Gamma(1-\alpha)}
    {}_{2}F_{1} \left[ 1,-n;1-\alpha; -\frac{x}{\beta} \right]
    \label{eq:RL derivative Analytical shifted poly}      
\end{split}
\end{equation}

\subsection*{Left-handed RL integral}

Due to the identity $\prescript{RL}{0}{\mathrm{D}}_x^{-\alpha}(\cdot) = \prescript{RL}{0}{\mathrm{I}}_x^{\alpha}(\cdot)$ $(0<\alpha<1)$ \cite{benkhettou2016existence}, the analytical solutions for the CO RL integral can be simply obtained from Eqs.~\ref{eq:RL derivative Analytical sin}-~\ref{eq:RL derivative Analytical shifted poly} by changing the sign of the order $\alpha$.

\clearpage
%\bibliographystyle{unsrt} 
%\bibliography{Ref}

\end{document}